# Output-only road roughness identification from vehicle axle accelerations through a universal smoothing method

Zihao Liu[a], Daniel Dias-da-Costa[b], Tommy Chan[c], Colin Caprani[d], Chul-Woo Kim[e], Mehrisadat Makki Alamdari[a,*]

[a]*School of Civil and Environmental Engineering, University of New South Wales, Sydney, NSW 2033, Australia*
[b]*School of Civil Engineering, The University of Sydney, Sydney, NSW 2006, Australia*
[c]*School of Civil & Environmental Engineering, Queensland University of Technology, Brisbane, QLD 4000, Australia*
[d]*Department of Civil and Environmental Engineering, Monash University, Clayton, VIC 3800, Australia*
[e]*Department of Civil & Earth Resources Engineering, Graduate School of Engineering, Kyoto University, Kyoto 615-8540, Japan*

**Abstract**

This paper presents an output-only method to identify road roughness profiles from axle accelerations of a moving vehicle. A two-degree-of-freedom half-car model is discretised with a zero-order hold and a backward-difference approximation of the roughness rate, which introduces both the current and previous roughness inputs into the observation equation. This modification enables joint input–state estimation with limited measurements using a Universal Smoothing (US) method, which belongs to the family of Minimum-Variance Unbiased (MVU) estimators. To improve numerical robustness under high process noise – stemming from modelling errors such as neglected bridge–vehicle interaction – the system inversion is regularised by truncated singular value decomposition. The method is validated on a full-scale bridge with a commercial SUV at two different speeds. Compared to the Dual Kalman filter and an MVU-based smoother, the proposed US achieves stable, accurate reconstructions across different scenarios and remains numerically well-conditioned when noise increases. Practical aspects of tuning, window length selection, and computational cost are also discussed.

*Corresponding author
*Email address:* m.makkialamdari@unsw.edu.au (Mehrisadat Makki Alamdari)

## Highlights

- A MVU-based recursive smoothing method is deployed for drive-by road roughness identification.
- A modified half-car model and its state-space representation are proposed to improve the invertibility and observability.
- The method is validated through a field test conducted at various car speeds.
- The proposed method demonstrates superior estimation quality over other state-of-the-art estimators despite high modelling error.

List of abbreviations

| | |
|---|---|
| AKF | Augmented Kalman filter |
| DKF | Dual Kalman filter |
| KF | Kalman filter |
| MVU | Minimum-variance unbiased |
| MVUS | Minimum-variance unbiased smoothing |
| NRMSE | Normalised root-mean-square error |
| UF | Universal filter |
| US | Universal smoothing |

## 1. Introduction

Road roughness (road profile) is a critical determinant of road-user safety and system performance: it influences rolling resistance, ride quality, vehicle operating costs, and fuel consumption, and it underpins asset management and vehicle/active-suspension design. To safeguard users and optimise maintenance, roughness must be located, measured, recorded, and characterised reliably—and done so in a timely, cost-effective manner [1]. The conventional road roughness estimation methods can be categorised into scanning-based and estimation-based methods. Imaging and scanning techniques have emerged as attractive, viable tools for inspecting and classifying pavement and road surface [2]. These methods can be contact-based or non-contact-based. Contact methods (e.g., rod-and-level, contact profilometers [3]) are inexpensive but labour-intensive, slow, and susceptible to observer-dependent errors, limiting scalability. Non-contact methods based on optical scanning mounted on survey vehicles (e.g., laser [1], LiDAR [4]) deliver high-precision profiles suitable for detailed mapping, yet their acquisition and maintenance costs are high, making routine or network-wide assessments on ordinary roads impractical. Estimation-based approaches [5, 6] infer the road profile indirectly from vehicle dynamics – typically using low-cost accelerometers or displacement sensors [7] – offering scalable, in-situ estimation [8]; however, inferring the road profile from vehicle responses is an ill-conditioned inverse problem and accuracy depends on vehicle modelling, parameter identification, operating conditions, and noise handling [9], motivating robust inverse-dynamics and system identification to improve reliability [10].

The Kalman Filter (KF) [11] is one of the most widely used recursive Bayesian filters in civil engineering applications, including road roughness identification. Early work applied the KF to a quarter-car with accelerometers and suspension-deflection sensors to reconstruct the road profile [12], establishing a practical baseline for estimation-based approaches. Augmented Kalman Filter (AKF) [13] enables joint input–state estimation by including the road profile as unknown input in the state vector. The AKF was then applied with a half-car model for road profile estimation using smartphone data, including acceleration, velocity, GPS and gyroscope [14]. Li et al. [15] utilised the AKF with a Bayesian expectation-maximisation method to reconstruct contact-point responses for indirect bridge roughness identification. A dual-wheeled scan-vehicle study was conducted by Li et al. [16] where bridge roughness is processed with a modified AKF to process limited data. Other refinements, e.g., trapezoidal load approximation with AKF, have also been proposed to enhance the road roughness identification [17]. However, a key limitation of the AKF is the need to assume a stochastic evolution for the unknown input within the augmented state vector. This assumption increases the sensitivity of the filter to process and measurement noise tuning, often forcing practitioners to rely on post-filter adjustments or richer sensor configurations. In practice, this frequently necessitates incorporating displacement-type measurements (e.g., suspension deflection) alongside acceleration data to maintain a well-conditioned estimation problem. To enable input estimation with acceleration-only measurement and mitigate the so-called drift issue, the DKF was developed by Eftekhar Azam et al. [18]. The DKF has been employed to reconstruct the road or bridge surface profile, which represents an essential preliminary step for bridge damage identification. Applications include a two-stage drive-by approach that first estimates surface roughness and subsequently detects structural damage [19]; a Newmark-beta-based DKF framework enabling simultaneous identification of road profiles and bridge support rotational stiffness [20], and a hybrid method that couples a Long Short-Term Memory network with the DKF for real-time, onboard road roughness estimation [21].

Motivated by the classical KF – which assumes known inputs – the Minimum-Variance Unbiased (MVU) filter [22] was developed to enable state estimation under unknown inputs. The original MVU, however, does not estimate the inputs themselves. To address this, two main variants have been proposed: (i) a formulation for systems without direct feedthrough [23] – in structural dynamics terms, this implies no acceleration measurements are used in the estimation; and (ii) a formulation for systems with direct feedthrough that requires full-rank feedforward [24], meaning the number of acceleration measurements must be at least equal to the number of unknown inputs. Both MVU estimators impose strict structural and rank conditions on the system and sensor layout. Even so, several studies have applied MVU methods to road roughness identification. For example, Shereena and Rao [25] used an MVU estimator to recover the road profile; to satisfy the system inversion and observability requirements, they combined vertical acceleration and deflection measurements and reported the need for a relatively long test section. Notably, their study is numerical and assumes known process and measurement noise, leaving real-world applicability uncertain. More recently, Kong et al. [26] developed an adaptive MVU framework based on a quarter-car model to estimate road roughness while simultaneously updating noise statistics and selected vehicle parameters during the estimation process. Despite these advances, MVU-based road roughness identification remains relatively unexplored in recent literature. Moreover, the rank and observability constraints of MVU frameworks impose strict requirements on sensor availability – often necessitating displacement/deflection measurements in addition to acceleration, particularly as the dimension of the unknown input increases.

A common limitation of the foregoing filtering methods is their reliance on measurements at a single timestep, which can render the estimation ill-conditioned when sensors are sparse, measurement noise is high, or model uncertainty is significant. A remedy is to employ smoothing methods that exploit measurements over an extended time window. By aggregating information forward and backward in time, smoothing effectively increases the available sensor information and regularises the inverse problem, yielding more accurate and robust state and input estimates. Maes et al. [27] introduced a smoothing counterpart to the MVU filter of [24] for systems with direct feedthrough – referred to hereafter as Minimum-Variance Unbiased Smoothing (MVUS) in the remainder of the paper. However, core MVU limitations remain, notably the requirement of a full-rank feedforward matrix. To overcome these constraints, the Universal Filter (UF) [28] and its generalisation, the Universal Smoothing (US) [29], were proposed. Prior work [30] suggests that UF/US can preserve well-conditioned system inversion across diverse sensor networks in structural health monitoring, yet their applicability to road-roughness identification has not been established.

This paper proposes an output-only method for identifying the road roughness profile using only axle accelerations based on the Universal Smoothing framework [29]. To make MVU-based estimation applicable under limited measurements, the temporal derivative of the roughness profile used in the half-car model is approximated by a backward difference, thereby reducing the dimension of the input vector and satisfying the minimum system inversion requirement of the MVU-based estimators. The modification introduces the previous-step roughness into both the process and observation equations, adding an additional source of uncertainty; accordingly, the weight matrix of the weighted least-squares involved in the US is augmented to propagate the covariance of the previous-step input estimate. Numerical robustness to high process noise arising from modelling errors, e.g., neglected bridge–vehicle interaction, is further improved by regularising the system inversion with truncated singular value decomposition. Validation on a full-scale bridge test – using a commercial SUV vehicle across two different speeds – benchmarks the approach against the DKF and MVUS, demonstrating stable, accurate reconstructions in cases where MVUS becomes ill-conditioned and DKF exhibits drift, and providing practical guidance on tuning and window-length selection.

This paper is organised as follows. Section 2 derives the two-degree-of-freedom half-car model, its discrete-time recursion, and introduces a backward-difference input that meets the invertibility conditions of MVU estimators under acceleration-only measurements; the implementation of the US is then described. Section 3 details the full-scale bridge field test, including site preparation, vehicle, instrumentation, and test scenarios. Section 4 reports road-roughness identification results and comparisons with DKF and MVUS, together with the tuning and regularisation settings used. Finally, conclusions and future research directions are provided in Section 5.

## 2. System Model and Smoothing Method

This section introduces the half-car model and the formulation of its dynamic system. The challenges regarding the road roughness identification through MVU estimators with limited measurements are discussed, and thereby, the modified state-space representation is introduced. Finally, the estimation steps and implementation of the US are explained.

### *2.1. Dynamic system of a road-vehicle model*

For the two-degree-of-freedom half-car model shown in Fig. 1, the formulation as in [31] is adopted herein. First, the equation of motion is given by,

$$\mathbf{M}_v\ddot{\mathbf{u}}(t) + \mathbf{C}_v\dot{\mathbf{u}}(t) + \mathbf{K}_v\mathbf{u}(t) = \mathbf{f}(t) \tag{1}$$

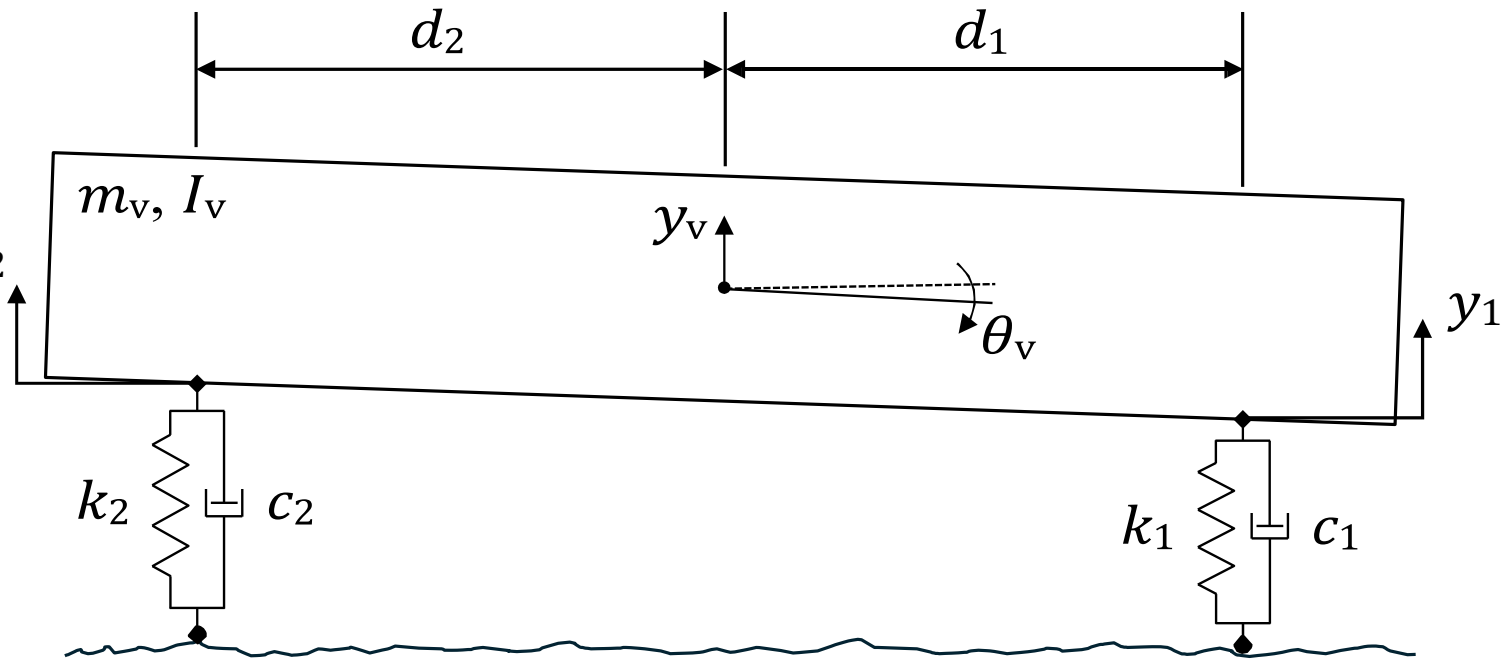

Figure 1: Two-degrees-of-freedom half-car model.

where the vector $\mathbf{u} \in \mathbb{R}^2$ contains the bounce and pitch displacements at the gravity centre. $\mathbf{M}_v \in \mathbb{R}^{2\times2}$, $\mathbf{C}_v \in \mathbb{R}^{2\times2}$ and $\mathbf{K}_v \in \mathbb{R}^{2\times2}$ are the mass, damping and stiffness matrices of the car, respectively. These vectors and matrices can be defined as follows,

$$\mathbf{u}(t) = [y_v(t)^T \quad \theta_v(t)^T]^T \tag{2}$$

$$\mathbf{M}_v = \begin{bmatrix} m_v & 0 \\ 0 & I_v \end{bmatrix} \tag{3}$$

$$\mathbf{K}_v = \begin{bmatrix} k_1 + k_2 & d_1k_1 - d_2k_2 \\ d_1k_1 - d_2k_2 & d_1^2k_1 + d_2^2k_2 \end{bmatrix}, \tag{4}$$

$$\mathbf{C}_v = \begin{bmatrix} c_1 + c_2 & d_1c_1 - d_2c_2 \\ d_1c_1 - d_2c_2 & d_1^2c_1 + d_2^2c_2 \end{bmatrix}. \tag{5}$$

In the above system matrices, $m_v$ is the mass and $I_v$ is the vehicle's moment of inertia, $k_1$, $c_1$, $k_2$, and $c_2$ are the stiffness and damping coefficients of the vehicle's front and rear axles, respectively. In addition, $d_1$ is the distance from the centre of gravity to the front axle, while $d_2$ is the distance to the rear axle. Furthermore, in Eq. (1), the external force vector $\mathbf{f}(t)$ is related to the road roughness profile and can be written as,

$$\mathbf{f}(t) = \mathbf{K}_r\mathbf{r}(t) + \mathbf{C}_r\dot{\mathbf{r}}(t), \tag{6}$$

in which, $\mathbf{r}$ is the road roughness profile vector that contains the roughness height under the front and rear axles. Note that the time function of the road roughness, $\mathbf{r}(t)$, can be obtained by using the velocity of the vehicle. Matrices $\mathbf{K}_r$ and $\mathbf{C}_r$ relate the vehicle properties to the road roughness, and read below,

$$\mathbf{K}_r = \begin{bmatrix} k_1 & k_2 \\ d_1k_1 & -d_2k_2 \end{bmatrix}, \tag{7}$$

$$\mathbf{C}_r = \begin{bmatrix} c_1 & c_2 \\ d_1c_1 & -d_2c_2 \end{bmatrix}. \tag{8}$$

The main challenge for deploying MVU-based approaches presented in the above half-car model is the system inversion requirement. In this study, the measured quantities are only the acceleration at the front and rear axles, i.e., $\ddot{y}_1$ and $\ddot{y}_2$. The system inversion requirement for the MVU filter with direct feedthrough [24] and its smoothing counterpart [27], MVUS, is that the number of measured accelerations is equal to or greater than the number of unknown inputs, while the requirement for the UF or US is that the number of any type of measurement is equal to or greater than the number of unknown inputs. However, the measured quantities in this study are the acceleration at the front and rear axles only. The limited measurement has violated the aforementioned minimum system inversion requirements. To solve this issue, Eq. (6) should be modified. First, the relation between $\mathbf{r}(t)$ and $\mathbf{r}(t)$ can be approximated by the backward differentiation,

$$\dot{\mathbf{r}}(t) = \frac{\mathbf{r}(t) - \mathbf{r}(t-\Delta t)}{\Delta t}. \tag{9}$$

Next, the second-order ordinary differential equation of Eq. (1) while considering Eq. (9) can be written as,

$$\dot{\mathbf{x}}(t) = \boldsymbol{A}\mathbf{x}(t) + \boldsymbol{B}\mathbf{r}(t) - \boldsymbol{G}\mathbf{r}(t-\Delta t), \tag{10}$$

where,

$$\mathbf{x}(t) = \begin{bmatrix} \mathbf{u}(t) \\ \dot{\mathbf{u}}(t) \end{bmatrix}, \tag{11}$$

$$\boldsymbol{A} = \begin{bmatrix} \mathbf{0}_{2\times 2} & \mathbf{I}_2 \\ -\mathbf{M}_v^{-1}\mathbf{K}_v & -\mathbf{M}_v^{-1}\mathbf{C}_v \end{bmatrix}, \tag{12}$$

$$\boldsymbol{B} = \begin{bmatrix} \mathbf{0}_{2\times 2} \\ \mathbf{M}_v^{-1}\mathbf{K}_r + \mathbf{M}_v^{-1}\mathbf{C}_r/\Delta t \end{bmatrix}, \tag{13}$$

$$\boldsymbol{G} = \begin{bmatrix} \mathbf{0}_{2\times 2} \\ \mathbf{M}_v^{-1}\mathbf{C}_r/\Delta t \end{bmatrix}, \tag{14}$$

and $\Delta t$ is the time interval. The time-discrete form of Eq. (10) can be obtained by employing the matrix exponential method [32],

$$\mathbf{x}_k = \mathbf{A}\mathbf{x}_{k-1} + \mathbf{B}\mathbf{r}_k - \mathbf{G}\mathbf{r}_{k-1}. \tag{15}$$

In the above equation, with a defined time interval $\Delta t$, the $k^{\text{th}}$ timestep can be represented as $t_k = k\Delta t$, thereby, $\mathbf{x}_k = \mathbf{x}(t_k)$, $\mathbf{r}_k = \mathbf{r}(t_k)$, $\mathbf{A} = \exp(\boldsymbol{A}\Delta t)$, $\mathbf{B} = \boldsymbol{B}\Delta t$ and $\mathbf{G} = \boldsymbol{G}\Delta t$.

### *2.2. State-space equations*

The state-space representation needs to be established in order to deploy system identification methods. Firstly, the process equation can be defined based on Eq. (15),

$$\mathbf{x}_k = \mathbf{A}_k\mathbf{x}_{k-1} + \mathbf{B}_{k-1}\mathbf{r}_k - \mathbf{G}_{k-1}\mathbf{r}_{k-1} + \mathbf{w}_{k-1}, \tag{16}$$

where $\mathbf{x}_k \in \mathbb{R}^{2n}$ is the state vector, $\mathbf{r}_k \in \mathbb{R}^m$ is the input vector; $\mathbf{A}_k \in \mathbb{R}^{2n\times 2n}$, $\mathbf{B}_k \in \mathbb{R}^{2n\times m}$ and $\mathbf{G}_k \in \mathbb{R}^{2n\times m}$ are the state and input matrices, respectively. $\mathbf{w}_k \in \mathbb{R}^{2n}$ stands for the process noise and is assumed to be zero-mean white noise. It should be noted that, for a linear system, the matrix $\mathbf{A}_k$, $\mathbf{B}_k$ and $\mathbf{G}_k$ are time-invariant in alignment with Eq. (15). However, their subscripts are kept in the smoothing formulation for the sake of completeness. Similarly, for the system size $n = 2$ and the input size $m = 2$, corresponding to the aforementioned two-degree-of-freedom half-car model, the symbols are used in the formulation for generality. Furthermore, the zero-order-hold adopted in the term $\mathbf{B}_{k-1}\mathbf{r}_k$ assumes that the input vector is expressed as $\mathbf{r}_k$, while $\mathbf{r}_{k-1}$ is commonly seen in literature [13, 23, 24], which leads to one timestep lag behind the input in the observation equation.

Next, assuming a quantity of $q$ sparse measurements on the system outputs are collected in the output vector $\mathbf{y}_k \in \mathbb{R}^q$, the observation equation can be expressed below,

$$\mathbf{y}_k = \mathbf{C}_k\mathbf{x}_k + \mathbf{D}_k\mathbf{r}_k - \mathbf{H}_k\mathbf{r}_{k-1} + \mathbf{v}_k, \tag{17}$$

in which the output matrix $\mathbf{C}_k \in \mathbb{R}^{q\times 2n}$, the feedforward matrices $\mathbf{D}_k \in \mathbb{R}^{q\times m}$ and $\mathbf{H}_k \in \mathbb{R}^{q\times m}$ relate the state vector $\mathbf{x}_k$, the input vector $\mathbf{r}_k$ and $\mathbf{r}_{k-1}$ to the output vector $\mathbf{y}_k$, respectively. The measurement noise $\mathbf{v}_k \in \mathbb{R}^q$ is also assumed to be zero-mean white noise. If the sensor network remains unchanged during the system identification process, the output matrix $\mathbf{C}_k$ is a constant and can be formulated as,

$$\mathbf{C}_k = \mathbb{C}\begin{bmatrix} \mathbf{I}_n & \mathbf{0}_{n\times n} \\ \mathbf{0}_{n\times n} & \mathbf{I}_n \\ -\mathbf{M}_v^{-1}\mathbf{K}_v & -\mathbf{M}_v^{-1}\mathbf{C}_v \end{bmatrix}, \tag{18}$$

where $\mathbb{C} \in \mathbb{R}^{q\times 3n}$ is a Boolean matrix containing the locations of sensors and the type of dynamic responses they collect. Similarly, the feedforward matrix $\mathbf{D}_k$, and the matrix $\mathbf{H}_k$ are also constant for the defined linear system,

$$\mathbf{D}_k = \mathbb{C}\begin{bmatrix} \mathbf{0}_{2n\times m} \\ \mathbf{M}_v^{-1}\mathbf{K}_r + \mathbf{M}_v^{-1}\mathbf{C}_r/\Delta t \end{bmatrix}. \tag{19}$$

$$\mathbf{H}_k = \mathbb{C}\begin{bmatrix} \mathbf{0}_{2n\times m} \\ \mathbf{M}_v^{-1}\mathbf{C}_r/\Delta t \end{bmatrix}. \tag{20}$$

The error covariance matrices of process and measurement noise can be defined as,

$$\begin{aligned} \mathbf{Q}_k &\triangleq \mathrm{E}\,[\mathbf{w}_k\mathbf{w}_k^T], \\ \mathbf{R}_k &\triangleq \mathrm{E}\,[\mathbf{v}_k\mathbf{v}_k^T]. \end{aligned} \tag{21}$$

Eqs. (16) and (17) can be adopted to most recursive filters. To deploy smoothing methods, an observation time window should be defined. Assuming a time window of length $N$, an extended observation equation can be formulated,

$$\boldsymbol{y}_k = \boldsymbol{C}_k\mathbf{x}_k + \boldsymbol{D}_k\boldsymbol{r}_k - \boldsymbol{H}_k\boldsymbol{r}_{k-1} + \boldsymbol{J}_k\boldsymbol{w}_k + \boldsymbol{v}_k, \tag{22}$$

where the extended vectors are defined as,

$$\begin{aligned} \boldsymbol{y}_k &\triangleq [\mathbf{y}_k^T \ \cdots \ \mathbf{y}_{k+N}^T]^T, \\ \boldsymbol{r}_k &\triangleq [\mathbf{r}_k^T \ \cdots \ \mathbf{r}_{k+N}^T]^T, \\ \boldsymbol{w}_k &\triangleq [\mathbf{w}_{k-1}^T \ \cdots \ \mathbf{w}_{k+N-1}^T]^T, \\ \boldsymbol{v}_k &\triangleq [\mathbf{v}_k^T \ \cdots \ \mathbf{v}_{k+N}^T]^T. \end{aligned} \tag{23}$$

and the extended matrices read,

$$\boldsymbol{C}_k = \begin{bmatrix} \mathbf{C}_k \\ \mathbf{C}_{k+1}\mathbf{A}_k \\ \mathbf{C}_{k+2}\mathbf{A}_{k+1}\mathbf{A}_k \\ \vdots \\ \mathbf{C}_{k+N}\prod_{i=0}^{N-1}\mathbf{A}_{k+i} \end{bmatrix}, \tag{24}$$

$$\boldsymbol{D}_k = \begin{bmatrix} \mathbf{D}_k & \mathbf{0} & \mathbf{0} & \cdots & \mathbf{0} \\ \mathbf{0} & \mathbf{C}_{k+1}\mathbf{B}_k + \mathbf{D}_{k+1} & \mathbf{0} & \cdots & \mathbf{0} \\ \mathbf{0} & \mathbf{C}_{k+2}\mathbf{A}_{k+1}\mathbf{B}_k & \mathbf{C}_{k+2}\mathbf{B}_{k+1} + \mathbf{D}_{k+2} & \cdots & \mathbf{0} \\ \mathbf{0} & \mathbf{C}_{k+3}\mathbf{A}_{k+2}\mathbf{A}_{k+1}\mathbf{B}_k & \mathbf{C}_{k+3}\mathbf{A}_{k+2}\mathbf{B}_{k+1} & \cdots & \mathbf{0} \\ \vdots & \vdots & \vdots & \ddots & \vdots \\ \mathbf{0} & \mathbf{C}_{k+N}\left(\prod_{i=1}^{N-1}\mathbf{A}_{k+i}\right)\mathbf{B}_k & \mathbf{C}_{k+N}\left(\prod_{i=2}^{N-1}\mathbf{A}_{k+i}\right)\mathbf{B}_{k+1} & \cdots & \mathbf{C}_{k+N}\mathbf{B}_{k+N-1} + \mathbf{D}_{k+N} \end{bmatrix}, \tag{25}$$

$$\boldsymbol{H}_k = \begin{bmatrix} \mathbf{H}_k & \mathbf{o} & \mathbf{o} & \cdots & \mathbf{o} \\ \mathbf{o} & \mathbf{C}_{k+1}\mathbf{G}_k + \mathbf{H}_{k+1} & \mathbf{o} & \cdots & \mathbf{o} \\ \mathbf{o} & \mathbf{C}_{k+2}\mathbf{A}_{k+1}\mathbf{G}_k & \mathbf{C}_{k+2}\mathbf{G}_{k+1} + \mathbf{H}_{k+2} & \cdots & \mathbf{o} \\ \mathbf{o} & \mathbf{C}_{k+3}\mathbf{A}_{k+2}\mathbf{A}_{k+1}\mathbf{G}_k & \mathbf{C}_{k+3}\mathbf{A}_{k+2}\mathbf{G}_{k+1} & \cdots & \mathbf{o} \\ \vdots & \vdots & \vdots & \ddots & \vdots \\ \mathbf{o} & \mathbf{C}_{k+N}\left(\prod_{i=1}^{N-1}\mathbf{A}_{k+i}\right)\mathbf{G}_k & \mathbf{C}_{k+N}\left(\prod_{i=2}^{N-1}\mathbf{A}_{k+i}\right)\mathbf{G}_{k+1} & \cdots & \mathbf{C}_{k+N}\mathbf{G}_{k+N-1} + \mathbf{H}_{k+N} \end{bmatrix}, \quad (26)$$

$$\boldsymbol{J}_k = \begin{bmatrix} \mathbf{0} & \mathbf{0} & \mathbf{0} & \cdots & \mathbf{0} \\ \mathbf{0} & \mathbf{C}_{k+1} & \mathbf{0} & \cdots & \mathbf{0} \\ \mathbf{0} & \mathbf{C}_{k+2}\mathbf{A}_{k+1} & \mathbf{C}_{k+2} & \cdots & \mathbf{0} \\ \mathbf{0} & \mathbf{C}_{k+3}\mathbf{A}_{k+2}\mathbf{A}_{k+1} & \mathbf{C}_{k+3}\mathbf{A}_{k+2} & \cdots & \mathbf{0} \\ \vdots & \vdots & \vdots & \ddots & \vdots \\ \mathbf{0} & \mathbf{C}_{k+N}\prod_{i=1}^{N-1}\mathbf{A}_{k+i} & \mathbf{C}_{k+N}\prod_{i=2}^{N-1}\mathbf{A}_{k+i} & \cdots & \mathbf{C}_{k+N} \end{bmatrix}. \quad (27)$$

Furthermore, the cross-covariance matrices of extended process and measurement noise can be defined as,

$$\begin{aligned} \boldsymbol{Q}_{i,j} &\triangleq \mathrm{E}\left[\boldsymbol{w}_i \boldsymbol{w}_j^T\right], \\ \boldsymbol{R}_{i,j} &\triangleq \mathrm{E}\left[\boldsymbol{v}_i \boldsymbol{v}_j^T\right]. \end{aligned} \quad (28)$$

*2.3. Smoothing steps*

The smoothing process of the US adopted is introduced in this section. At the start of each timestep, a biased state estimation is calculated given that the input is unknown,

$$\hat{\mathbf{x}}_{k|k-1} = \mathbf{A}_{k-1}\hat{\mathbf{x}}_{k-1|k-1} - \mathbf{G}_{k-1}\mathbf{r}_{k-1}, \quad (29)$$

In the equation above, the bias is caused by the missing contribution of the input. Therefore, the input can be inversely obtained by the innovation and a gain matrix $\mathbf{M}_k$. Note that the innovation in the smoothing method is calculated using the extended output vector,

$$\hat{\boldsymbol{r}}_k = \mathbf{M}_k\left(\boldsymbol{y}_k - \boldsymbol{C}_k\hat{\mathbf{x}}_{k|k-1} + \boldsymbol{H}_k\hat{\boldsymbol{r}}_{k-1}\right), \quad (30)$$

and the input vector $\hat{\mathbf{r}}_k$ can be extracted from $\hat{\boldsymbol{r}}_k$ using,

$$\hat{\mathbf{r}}_k = \left[\mathbf{I}_m \quad \mathbf{0}_{m\times Nm}\right]\hat{\boldsymbol{r}}_k. \quad (31)$$

Next, a priori state estimation can be obtained by adding the contribution of the estimated input,

$$\hat{\boldsymbol{x}}_{k|k-1} = \hat{\mathbf{x}}_{k|k-1} + \mathbf{B}_{k-1}\hat{\mathbf{r}}_k. \quad (32)$$

Finally, a posteriori state estimation can be obtained from the measurement in the time window with a state gain matrix $\mathbf{K}_k$,

$$\hat{\mathbf{x}}_{k|k} = \hat{\boldsymbol{x}}_{k|k-1} + \mathbf{K}_k\left(\boldsymbol{y}_k - \boldsymbol{C}_k\hat{\boldsymbol{x}}_{k|k-1} - \boldsymbol{D}_k\hat{\boldsymbol{r}}_k + \boldsymbol{H}_k\hat{\boldsymbol{r}}_{k-1}\right). \quad (33)$$

It should be noted that in the original formulation of the US, the weight matrix $\tilde{\mathbf{R}}_k$ for the weighted least squares is,

$$\tilde{\mathbf{R}}_k \triangleq \mathrm{E}\left[\mathbf{e}_k \mathbf{e}_k^T\right] = \boldsymbol{\Sigma}_k \boldsymbol{\Lambda}_k \boldsymbol{\Sigma}_k^T, \tag{34}$$

in which $\mathbf{e}_k$ is an error term that includes estimation error from state, input, process noise and measurement noise; matrix $\boldsymbol{\Lambda}_k$ contains all error covariance matrices involved in the estimation process,

$$\boldsymbol{\Lambda}_k \triangleq \begin{bmatrix} \mathbf{P}_{k-1} & \mathbf{P}^{\mathbf{xw}}_{k-1} & \mathbf{P}^{\mathbf{xv}}_{k-1} \\ \mathbf{P}^{\mathbf{wx}}_{k-1} & \boldsymbol{Q}_{k,k} & \mathbf{P}^{\mathbf{wv}}_{k} \\ \mathbf{P}^{\mathbf{vx}}_{k-1} & \mathbf{P}^{\mathbf{vw}}_{k} & \boldsymbol{R}_{k,k} \end{bmatrix}, \tag{35}$$

and matrix $\boldsymbol{\Sigma}_k$ contains system matrices linked to error covariance matrices,

$$\boldsymbol{\Sigma}_k \triangleq \left[\boldsymbol{\Gamma}_k \quad \breve{\mathrm{H}}_k \quad \mathbf{I}_{(N+1)q}\right], \tag{36}$$

where $\boldsymbol{\Gamma}_k = \boldsymbol{C}_k \mathbf{A}_{k-1}$ and $\breve{\mathrm{H}}_k \triangleq \boldsymbol{J}_k \left[+ \;\; \boldsymbol{C}_k \;\; \mathbf{0}_{(N+1)q\times N\times 2n}\right]$. The detailed derivation and the proof of unbiasedness can be found in [29]. Since the state-space representation established in this study introduces the input from the previous timestep, i.e., $\mathbf{r}_{k-1}$, to account for estimation error presented in $\mathbf{r}_{k-1}$, a new weight matrix is proposed,

$$\tilde{\mathbf{R}}_k \triangleq \mathrm{E}\left[\mathbf{e}_k \mathbf{e}_k^T\right] = \boldsymbol{\Sigma}_k \boldsymbol{\Lambda}_k \boldsymbol{\Sigma}_k^T + \boldsymbol{\Xi}_k \mathbf{P}^{\mathbf{r}}_{k-1} \boldsymbol{\Xi}_k^T, \tag{37}$$

where,

$$\boldsymbol{\Xi}_k = \begin{bmatrix} \mathbf{C}_k \mathbf{G}_{k-1} + \mathbf{H}_k \\ \mathbf{C}_{k+1} \mathbf{A}_k \mathbf{G}_{k-1} \\ \mathbf{C}_{k+2} \mathbf{A}_{k+1} \mathbf{A}_k \mathbf{G}_{k-1} \\ \vdots \\ \mathbf{C}_{k+N} \left(\prod_{i=1}^{N-1} \mathbf{A}_{k+i}\right) \mathbf{G}_{k-1} \end{bmatrix} \tag{38}$$

The error propagation and calculation of gain matrices are summarised in Table 1 for the ease of implementation. Note that $(\cdot)^{\dagger}$ represents the Moore–Penrose pseudoinverse.

Table 1: Summary of the Universal Smoothing

---

**Initialisation**

- Define process and observation equations
- Assign initial values to $\hat{\mathbf{x}}_{0|0}$, $\mathbf{P}_0$, $\mathbf{Q}_0$, $\mathbf{R}_0$ and build $\boldsymbol{Q}_{0,0}$, $\boldsymbol{Q}_{0,1}$, $\boldsymbol{R}_{0,0}$ and $\boldsymbol{R}_{0,1}$
- Assign $\mathbf{P}^{\mathbf{xw}}_0 = \mathbf{O}$ and $\mathbf{P}^{\mathbf{xv}}_0 = \mathbf{O}$
- Build time-independent matrices

$$\breve{\boldsymbol{\varepsilon}}_i \triangleq \begin{bmatrix} \mathbf{O}_{iN\times i} & \mathbf{I}_{Ni} \\ \mathbf{O}_{i\times i} & \mathbf{O}_{i\times iN} \end{bmatrix},\ \bar{\boldsymbol{\varepsilon}}_i \triangleq \begin{bmatrix} \mathbf{O}_{iN\times iN} & \mathbf{O}_{iN\times i} \\ \mathbf{O}_{i\times iN} & \mathbf{I}_i \end{bmatrix} \text{ for } i = 2n \text{ and } q$$

---

**Smoothing loop**

- Input estimation
  1. $\mathbf{R}_k = \mathbf{\Sigma}_k\mathbf{\Lambda}_k\mathbf{\Sigma}_k^T + \mathbf{\Xi}_k\mathbf{P}^{\mathbf{r}}_{k-1}\mathbf{\Xi}_k^T$
  2. $\mathbf{P}^{\boldsymbol{r}}_k = (\breve{\mathbf{D}}_k^T\tilde{\mathbf{R}}_k^{\dagger}\breve{\mathbf{D}}_k)^{\dagger}$
  3. $\mathbf{P}^{\mathbf{r}}_k = [\mathbf{I}_m \quad \mathbf{O}_{m\times Nm}]\,\mathbf{P}^{\boldsymbol{r}}_k\,[\mathbf{I}_m \quad \mathbf{O}_{m\times Nm}]^T$
  4. $\mathbf{M}_k = \mathbf{P}^{\boldsymbol{r}}_k\breve{\mathbf{D}}_k^T\tilde{\mathbf{R}}_k^{\dagger}$
  5. $\hat{\mathbf{x}}_{k|k-1} = \mathbf{A}_{k-1}\hat{\mathbf{x}}_{k-1|k-1} - \mathbf{G}_{k-1}\mathbf{r}_{k-1}$
  6. $\hat{\boldsymbol{r}}_k = \mathbf{M}_k(\boldsymbol{y}_k - \boldsymbol{C}_k\hat{\mathbf{x}}_{k|k-1} + \boldsymbol{H}_k\hat{\boldsymbol{r}}_{k-1})$
  7. $\hat{\mathbf{r}}_k = [\mathbf{I}_m \quad \mathbf{O}_{m\times Nm}]\,\hat{\boldsymbol{r}}_k$
- State estimation
  1. $\hat{\boldsymbol{x}}_{k|k-1} = \hat{\mathbf{x}}_{k|k-1} + \mathbf{B}_{k-1}\hat{\mathbf{r}}_k$
  2. $\mathbf{V}_k \triangleq \mathbf{B}_{k-1}[\mathbf{I}_m \quad \mathbf{O}_{m\times Nm}]\,\mathbf{M}_k$
  3. $\mathbf{W}_k \triangleq -\mathbf{V}_k\breve{\mathbf{H}}_k + [\mathbf{I}_{2n} \quad \mathbf{O}_{2n\times N\times 2n}]$
  4. $\breve{\mathbf{A}}_{k-1} \triangleq \mathbf{A}_{k-1} - \mathbf{V}_k\mathbf{\Gamma}_k$
  5. $\mathbf{\Pi}_k \triangleq [\breve{\mathbf{A}}_{k-1} \quad \mathbf{W}_k \quad -\mathbf{V}_k]$
  6. $\mathbf{P}^{\boldsymbol{x}}_k = \mathbf{\Pi}_k\mathbf{\Lambda}_k\mathbf{\Pi}_k^T$
  7. $\mathbf{\Theta}_k = \boldsymbol{D}_k\mathbf{M}_k$
  8. $\mathbf{\Omega}_k \triangleq [\boldsymbol{C}_k\breve{\mathbf{A}}_{k-1} - \mathbf{\Theta}_k\mathbf{\Gamma}_k \quad \boldsymbol{C}_k\mathbf{W}_k - \mathbf{\Theta}_k\breve{\mathbf{H}}_k + \boldsymbol{J}_k \quad -\boldsymbol{C}_k\mathbf{V}_k - \mathbf{\Theta}_k + \mathbf{I}_{(N+1)q\times 2n}]$
  9. $\mathbf{Y}_k = -\mathbf{\Omega}_k\mathbf{\Lambda}_k\mathbf{\Pi}_k^T$
  10. $\mathbf{\Phi}_k = \mathbf{\Omega}_k\mathbf{\Lambda}_k\mathbf{\Omega}_k^T$
  11. $\mathbf{K}_k = -\mathbf{Y}_k^T\mathbf{U}_k^T(\mathbf{U}_k\mathbf{\Phi}_k\mathbf{U}_k^T)^{-1}\mathbf{U}_k$
  12. $\hat{\mathbf{x}}_{k|k} = \hat{\boldsymbol{x}}_{k|k-1} + \mathbf{K}_k(\boldsymbol{y}_k - \boldsymbol{C}_k\hat{\boldsymbol{x}}_{k|k-1} - \boldsymbol{D}_k\hat{\boldsymbol{r}}_k + \boldsymbol{H}_k\hat{\boldsymbol{r}}_{k-1})$
  13. $\mathbf{P}_k = \mathbf{P}^{\boldsymbol{x}}_k + \mathbf{K}_k\mathbf{Y}_k + \mathbf{Y}_k^T\mathbf{K}_k^T + \mathbf{K}_k\mathbf{\Phi}_k\mathbf{K}_k^T$
  14. $\mathbf{T}_k = \mathbf{I}_{(N+1)q\times 2n} - \mathbf{K}_k\boldsymbol{C}_k$
  15. $\boldsymbol{W}_k = \mathbf{T}_k\mathbf{W}_k + \mathbf{K}_k\mathbf{\Theta}_k\breve{\mathbf{H}}_k - \mathbf{K}_k\boldsymbol{J}_k$
  16. $\boldsymbol{V}_k = -\mathbf{T}_k\mathbf{V}_k + \mathbf{K}_k\mathbf{\Theta}_k - \mathbf{K}_k$
  17. $\tilde{\boldsymbol{A}}_k = \mathbf{T}_k\breve{\mathbf{A}}_{k-1} + \mathbf{K}_k\mathbf{\Theta}_k\mathbf{\Gamma}_k$
  18. $\mathbf{P}^{\mathbf{xw}}_k = \tilde{\boldsymbol{A}}_k\mathbf{P}^{\mathbf{xw}}_{k-1}\breve{\boldsymbol{\varepsilon}}_{2n}^T + \boldsymbol{W}_k\boldsymbol{Q}_{k,k+1}$
  19. $\mathbf{P}^{\mathbf{xv}}_k = \tilde{\boldsymbol{A}}_k\mathbf{P}^{\mathbf{xv}}_{k-1}\breve{\boldsymbol{\varepsilon}}_q^T + \boldsymbol{V}_k\boldsymbol{R}_{k,k+1}$

---

## 3. Case Study: Road Roughness Identification Using Drive-by Inspection

The field experiment conducted by Chang et al. [31] and Hasegawa et al. [33] is used to validate the presented road roughness identification method. A passenger vehicle traversed the bridge multiple times at a nominally constant speed, with axle accelerations recorded during each pass. The road roughness profile was measured independently using a laser-based instrument as a reference (ground truth) to assess the estimation accuracy. Key aspects of the test site, vehicle, instrumentation, test procedure, data acquisition, and reference road profile are detailed below.

### *3.1. Test preparation: bridge, vehicle, instrumentation and road roughness profiles*

A drive-by test was conducted on a single-lane bridge with a span of approximately 40.5 m, located in Uji City, Kyoto. Fig. 2a shows the site map where the road profile beneath each side of the wheel is referred to as Track 1 and Track 2, respectively. The tracks were marked with red paint to ensure the wheels remained on the same track throughout multiple travels, as shown in Fig. 2b. The laser-measured road roughness profiles for each track are shown in Fig. 3 (measured from north to south). The vehicle was a Nissan X-Trail, modelled as the half-car model, as introduced in Section 2, for analysis. Representative parameters used in the model are summarised in Table 2. The properties were assembled from manufacturer information as well as free-vibration testing, and were reported in [31, 33]. In this case study, given the bridge dimensions, the bridge mass is assumed to be much larger than the vehicle mass. The vehicle-bridge interaction is therefore neglected, and the half-car model in Section 2 is used to relate axle accelerations to the road roughness. The measured accelerations still contain vibration from the bridge; this mismatch is treated as modelling error (process noise) in the estimators. This simplification deviates from the fully coupled problem, but it reflects common field conditions in which bridge properties and boundary conditions are unknown or cannot be determined in advance. It therefore imposes a demanding robustness requirement – estimators must operate under comparatively large process noise – while enabling a practical deployment setting in which a vehicle-only model and acceleration-only measurements can be used in general practice.

Table 2: Properties of the inspection vehicle.

| Properties | Values |
|---|---|
| Front wheel to centroid distance ($d_1$) | 0.82 m |
| Rear wheel to centroid distance ($d_2$) | 1.90 m |
| Mass ($m_v$) | 1, 994 kg |
| Moment of inertia ($I_v$) | 3, 954 kg · m$^2$ |
| Front axle spring constant $k_1$ | 75, 749 N/m |
| Rear axle spring constant $k_2$ | 99, 646 N/m |
| Front axle damping coefficient $c_1$ | 12, 535 N · s/m |
| Rear axle damping coefficient $c_2$ | 3, 602 N · s/m |

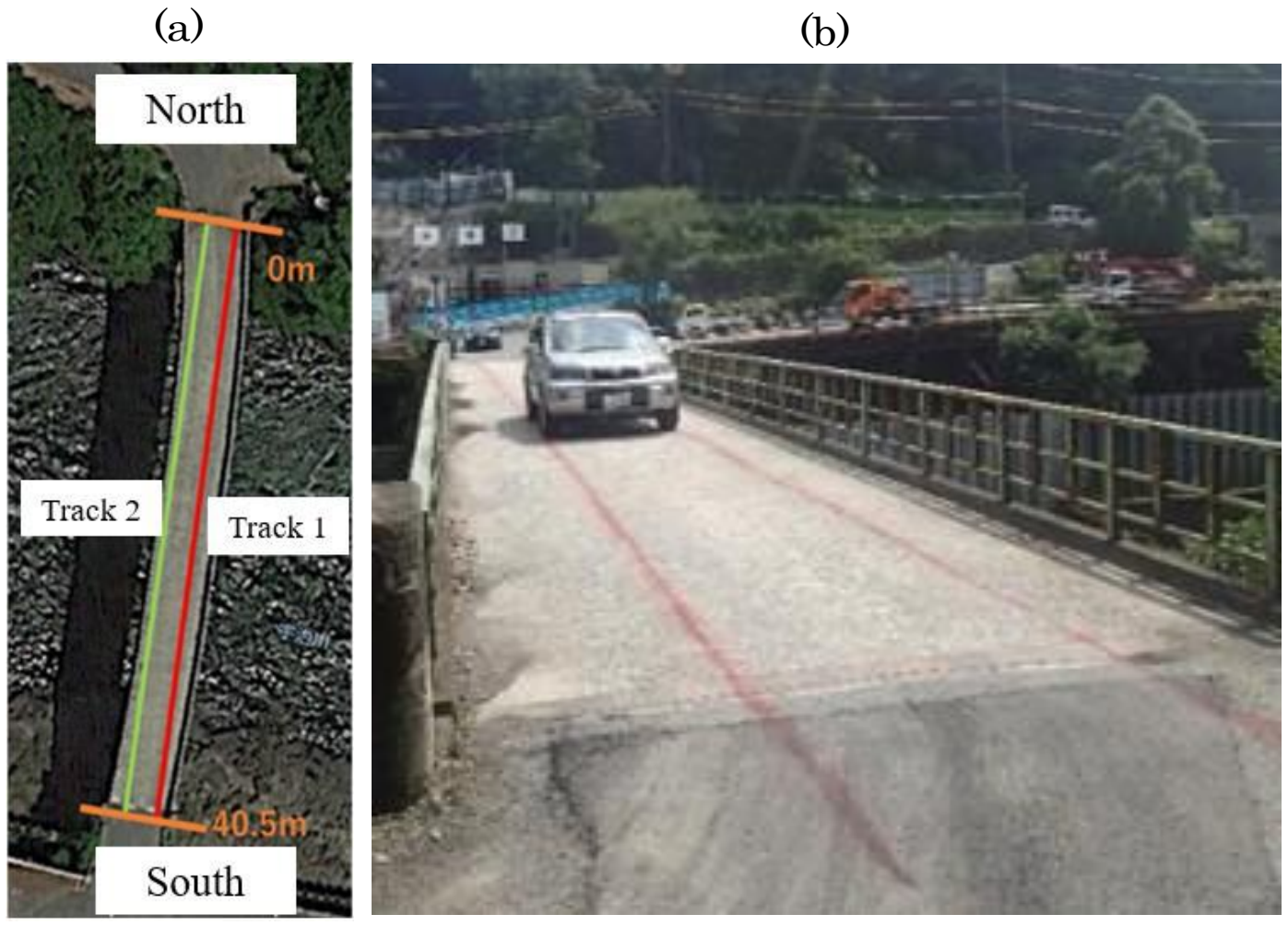

Figure 2: (a) Site map of the tested bridge and labels of road roughness profiles [31], and (b) the testing vehicle travelling on the marked tracks, (Photo credit: Kyoto University).

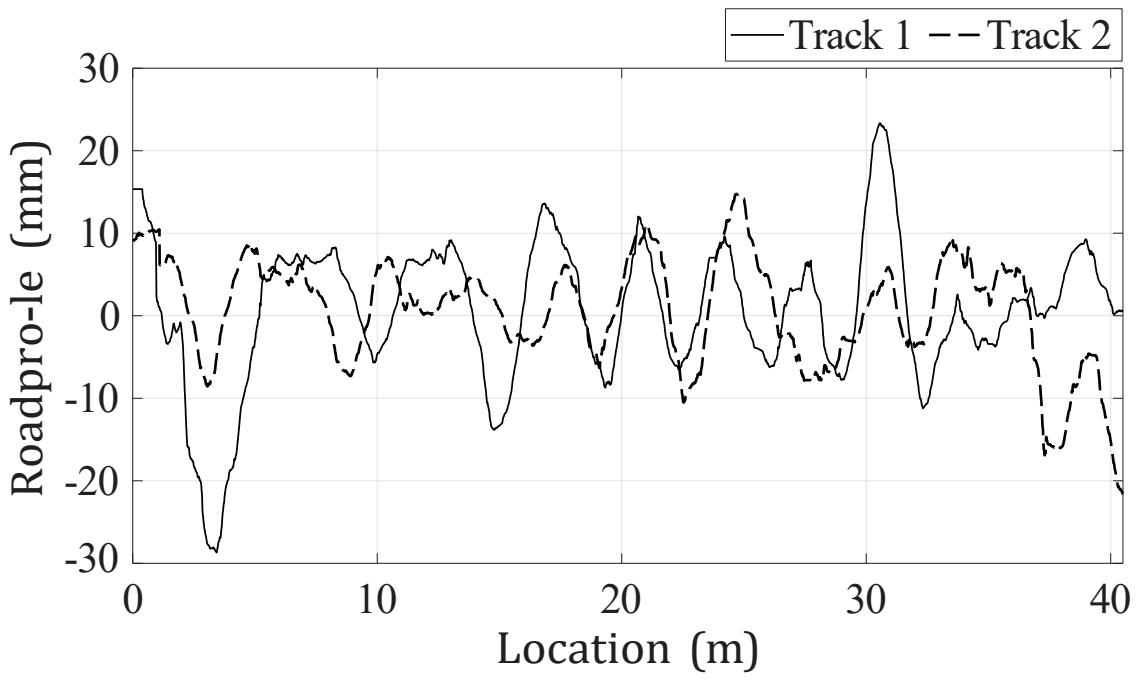

Figure 3: The laser measured road roughness profiles (from north to south).

Four accelerometers, Epson M-A550AC2, were mounted near the front and rear axles of the car to record the vertical acceleration with a sampling rate of 200 Hz. To suppress low-frequency drift and high-frequency noise while preserving the physically relevant vehicle–roughness content, the acceleration time series were band-pass filtered with cut-off frequencies of 0.1-50 Hz. The vertical accelerations obtained were processed using the rigid body assumption to determine the accelerations in the bounce and pitch directions required by the half-car model. It should be noted that, since the accelerometers did not exactly align with the axles, the actual distance to the centroid was used in the calculation, as shown in Fig 4. With the given sensor layout, if the vehicle were driven northbound, sensors 1 and 2 could be used to identify track 2, while track 1 could be identified by sensors 3 and 4.

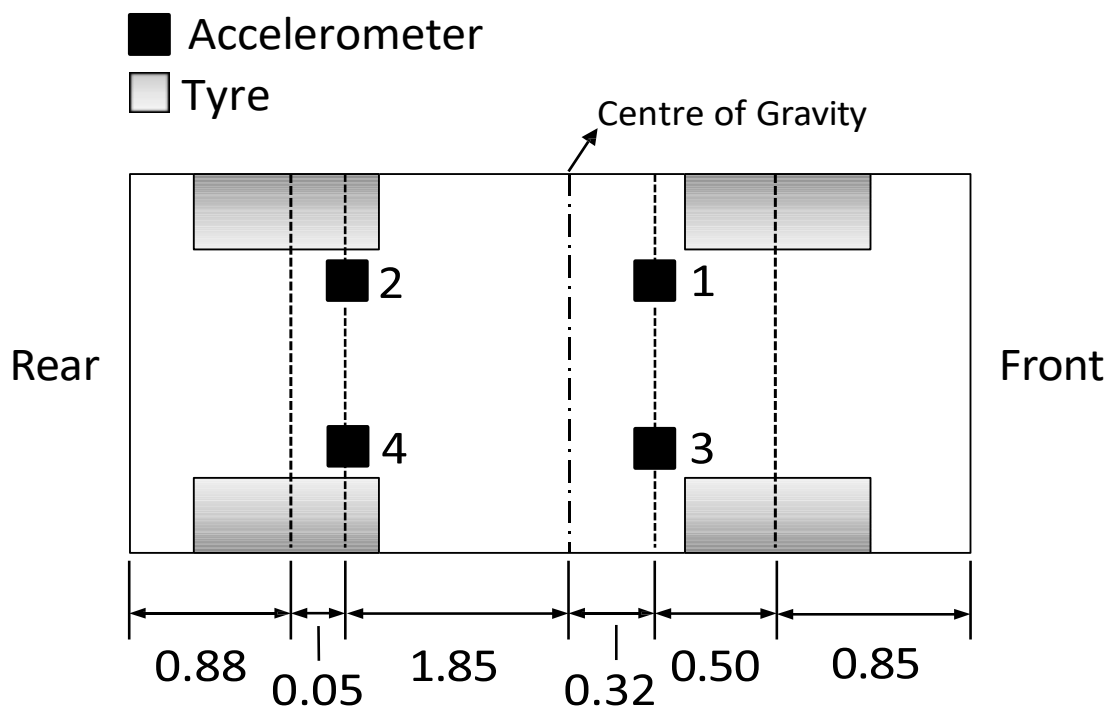

Figure 4: Sensor layout (unit: m, not to scale).

### 3.2. Identification preparation: scenarios, tuning strategies, and performance evaluation

Two different vehicle speeds of 10 km/h and 20 km/h were considered in the tests to examine the influence of the vehicle speed on the road roughness identification results. Multiple test runs were conducted during the experiment. The selected datasets and test scenarios are summarised in Table 3.

Table 3: Test scenarios used for road roughness profile identification.

| Scenario No. | Identified profile | Vehicle speed |
|---|---|---|
| 1 | Track 1 | 10 km/h |
| 2 | Track 1 | 20 km/h |
| 3 | Track 2 | 10 km/h |
| 4 | Track 2 | 20 km/h |

For recursive estimators, the measurement noise matrix $\mathbf{R}_k$, the process noise matrix $\mathbf{Q}_k$ and several tuning parameters should be set properly to achieve optimal results. The diagonal values of the matrix $\mathbf{R}_k$ are set based on the signal-to-noise ratio of the sensors. Since real-time estimation is not a concern in road roughness identification, an offline grid search is adopted to tune parameters. To reflect the practical condition where the true road roughness is unknown and only axle acceleration measurements are available, the estimation error proposed by Vettori et al. [34] is used in the tuning process, which is defined as,

$$\Sigma \mathbf{E} = \left\| [\mathbf{E}_u \quad \mathbf{E}_r]^T \right\|_{\frac{1}{2}} \tag{39}$$

In the above equation, $\mathbf{E}_u$ and $\mathbf{E}_r$ are the errors of estimated acceleration and road roughness, respectively. Given that accelerations are measured, $\mathbf{E}_u$ can be formulated as:

$$\mathbf{E}_u = \sqrt{\sum_{k=1}^{T} \left( \arg\min_{\beta} \|\beta \mathbf{u}_k - (\mathbf{u}_k - \hat{\mathbf{u}}_k)\|_2 \right)^2}, \tag{40}$$

where $T$ is the total timesteps, $\mathbf{u}_k$ is the measured acceleration, $\hat{\mathbf{u}}_k$ is the estimated acceleration, and $\beta$ is the scalar result of the least-squares solution of the equation $\beta u_k = u_k - \hat{u}_k$. This regulation ensures that the error term $E_u$ is dimensionless. The input is assumed to be unknown during the estimation process. Therefore, the error can be estimated based on the error covariance matrix $\mathbf{P}^{\mathbf{r}}_k$ and the squared value of the estimated input,

$$\mathbf{E}_r = \sqrt{\sum_{k=1}^{T} \frac{\mathrm{diag}(\mathbf{P}^{\mathbf{r}}_k)}{\mathbf{r}^2_k}}, \tag{41}$$

in which, $\mathrm{diag}(\mathbf{P}^{\mathbf{r}}_k)$ are the diagonal elements of the error covariance $\mathbf{P}^{\mathbf{r}}_k$, and the squared value $\mathbf{r}^2_k$ is used for normalisation. The tuning procedures for the US and the other methods used for comparison are discussed in Section 4.1. Finally, to accurately evaluate the performance across different methods, the Normalised Root Mean Square Error (NRMSE) is used to calculate the true estimation error (rather than estimated by $\Sigma\mathbf{E}$),

$$\mathrm{NRMSE} = \sum_{k=1}^{T} \left[ \frac{\Delta r_{rms,k}}{(r_{max} - r_{min})} \right], \tag{42}$$

where $r_{rms,k}$ is the root-mean-square error between the estimated and the true roughness at each timestep; $r_{max}$ and $r_{min}$ are the absolute maximum and minimum true values of road roughness, respectively.

## 4. Road Roughness Identification Results

This section presents road roughness identification results from the drive-by experiment. We benchmark the proposed US against two state-of-the-art methods: (i) the DKF [18], which mitigates drift under acceleration-only measurements and is noted for numerical stability, and (ii) the MVUS [27], a pioneering smoothing method within the MVU family. Section 4.1 outlines the tuning and initialisation strategy to ensure a fair comparison. The remaining sections report results for varying speed, an NRMSE-based quantitative comparison, and a discussion of estimator-specific behaviour.

### 4.1. Tuning, regularisation and initialisation

To tune the process noise $\mathbf{Q}_k$, all diagonal components are assumed identical; thus a scalar $Q_x$ is used such that $\mathbf{Q}_k = Q_x\mathbf{I}_4$. The scalar $Q_x$ is searched logarithmically from $10^{-12}$ to $10^{-1}$ with a step of $10^{0.1}$. For the DKF, an additional error matrix for the input vector is needed: $\mathbf{Q}_r = Q_r\,\mathbf{I}_2$, and $Q_r$ is searched from $10^{-8}$ to $10^{2}$ with a step ratio of $10^{0.1}$.

Compared to other applications, such as structural health monitoring, where the observed system is well characterised, the present study adopts a simplified half-car model and does not include bridge–vehicle interaction; i.e., excitation arising from bridge vibration is not modelled. Consequently, process noise representing modelling errors and uncertainties is significant. With limited measurements, numerical stability becomes challenging, particularly for MVU-based approaches. One advantage of the US is that the system inversion uses the Moore–Penrose pseudoinverse, enabling numerical stabilisation through singular-value truncation [35]. Truncation can be controlled either by a tolerance or by retaining a fixed number of singular values; the latter is adopted here, with the optimal count searched from 1 up to $2(N+1)$, corresponding to the maximum rank of the extended input vector $\hat{\boldsymbol{r}}_k$ across a window of length $N$.

Furthermore, the initialisation also affects the estimation quality. In all comparative studies presented later the initial state $\mathbf{x}_0$ is set to zero; the initial state error covariance $\mathbf{P}_0^{\mathbf{x}}$ is set to $10^{-12}\mathbf{I}_4$ and the input error covariance $\mathbf{P}_0^{\mathbf{r}}$ is set to $10^{-12}\mathbf{I}_2$. For smoothing methods—including the US and the MVUS—the time-window length N strongly influences performance by determining how much information the estimator can access at each timestep. However, as N increases, the dimensions of the extended system matrices grow rapidly, incurring a substantial computational cost. The impact of time window size $N$ on estimation accuracy and computational cost is presented in Fig. 5. For reference, the computation time is measured based on an AMD Ryzen 7 9800X3D Processor with 64 GB of memory. As shown in Fig. 5, improvement in estimation accuracy becomes marginal when $N > 100$, whereas the computational cost continues to increase exponentially. To balance accuracy and computation time, a window size of $N = 100$ is adopted for all smoothing methods.

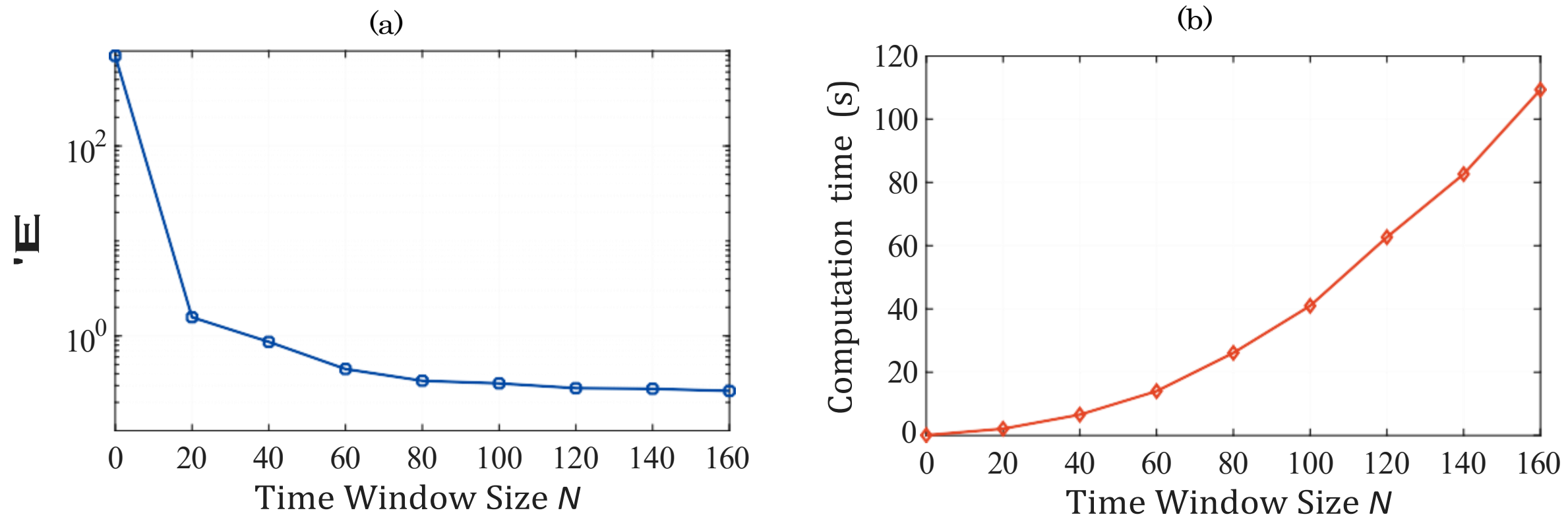

Figure 5: Impact of the time window size $N$ on the: (a) estimation quality, and (b) computation time.

With all settings discussed above, grid searches are conducted over the process noise scale and the truncation level used in the Moore–Penrose pseudoinverse (number of singular values retained). Tuning is performed independently for each scenario. The selected grid-search results for the US in Scenarios 2 and 4 are shown in Fig. 6. In both cases, a broad low-error region appears for small process-noise levels combined with moderate truncation. When only a few singular values are retained, the US remains numerically stable throughout the scanned $Q_x$ range and $\Sigma\mathbf{E}$ varies continuously. However, retaining only a small set of singular values removes the weakly observable components of the signal. This suppresses details and introduces a bias, resulting in high estimation error, even though the solution stays numerically stable. Conversely, as the number of retained singular values increases, the inverse becomes sensitive to process noise and $\Sigma\mathbf{E}$ exhibits sharp rises around $Q_x \gtrsim 10^{-6}$ before encountering numerical error. The optimal combination of tuned parameters is typically found just before ill-conditioning occurs. The tuning results are summarised in Table 4.

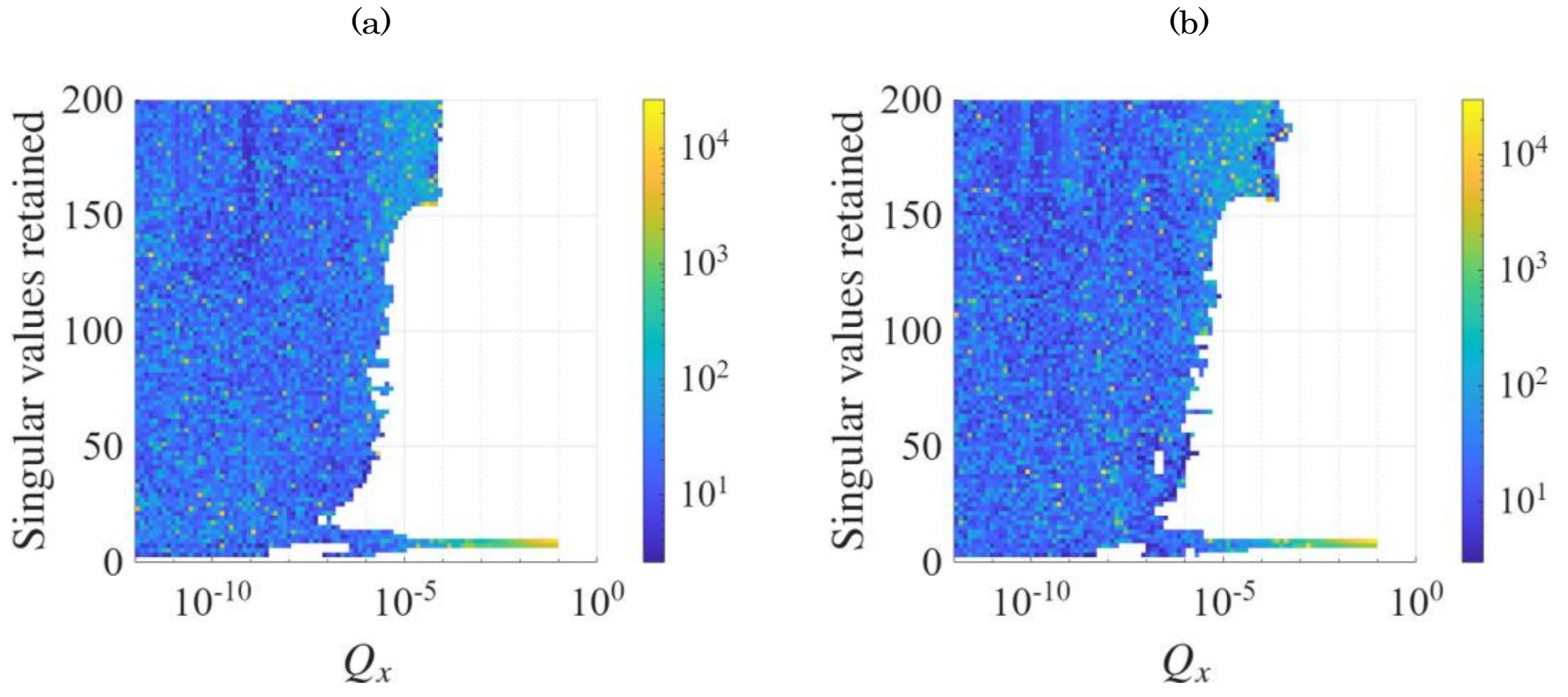

Figure 6: Selected grid search results of the US: (a) Scenario 2 and (b) Scenario 4.

Table 4: Optimal parameters for the US within the searched range.

| Scenario No. | $Q_x$ | Singular values retained | $\Sigma\mathbf{E}$ |
|---|---|---|---|
| 1 | $10^{-6.6}$ | 44 | 0.216 |
| 2 | $10^{-8.2}$ | 40 | 0.277 |
| 3 | $10^{-6.4}$ | 46 | 0.259 |
| 4 | $10^{-8.7}$ | 34 | 0.296 |

### *4.2. Identification results*

Fig. 7 presents the Track 1 profile identified at approximately 10 km/h (Scenario 1). The compared methods reproduce the overall roughness trend but uniformly underestimate the dominant peak occurring between 25 m and 30 m. In contrast, the DKF exhibits a more pronounced underestimation, aligning with the anticipated distinction between filters and smoothers: smoothing exploits measurements within a time window, providing more information per timestep. The corresponding NRMSE values are 0.216 for the US, 0.252 for the MVUS, and 0.304 for the DKF.

It is expected that a higher speed will intensify the measurement noise. To assess this effect, the Track 1 profile identified at 20 km/h is presented in Fig. 8 (Scenario 2). As the noise level increases, the grid search fails to identify any process-noise value within the prescribed range that allows the MVUS to proceed without a singular matrix inversion. This sensitivity to measurement noise and the associated numerical instability of MVU-type methods has also been reported in prior studies [28, 29, 30]. The DKF likewise begins to drift, whereas the US remains stable. The resulting NRMSE for the US is 0.218, indicating that its accuracy is largely preserved despite the higher noise level.

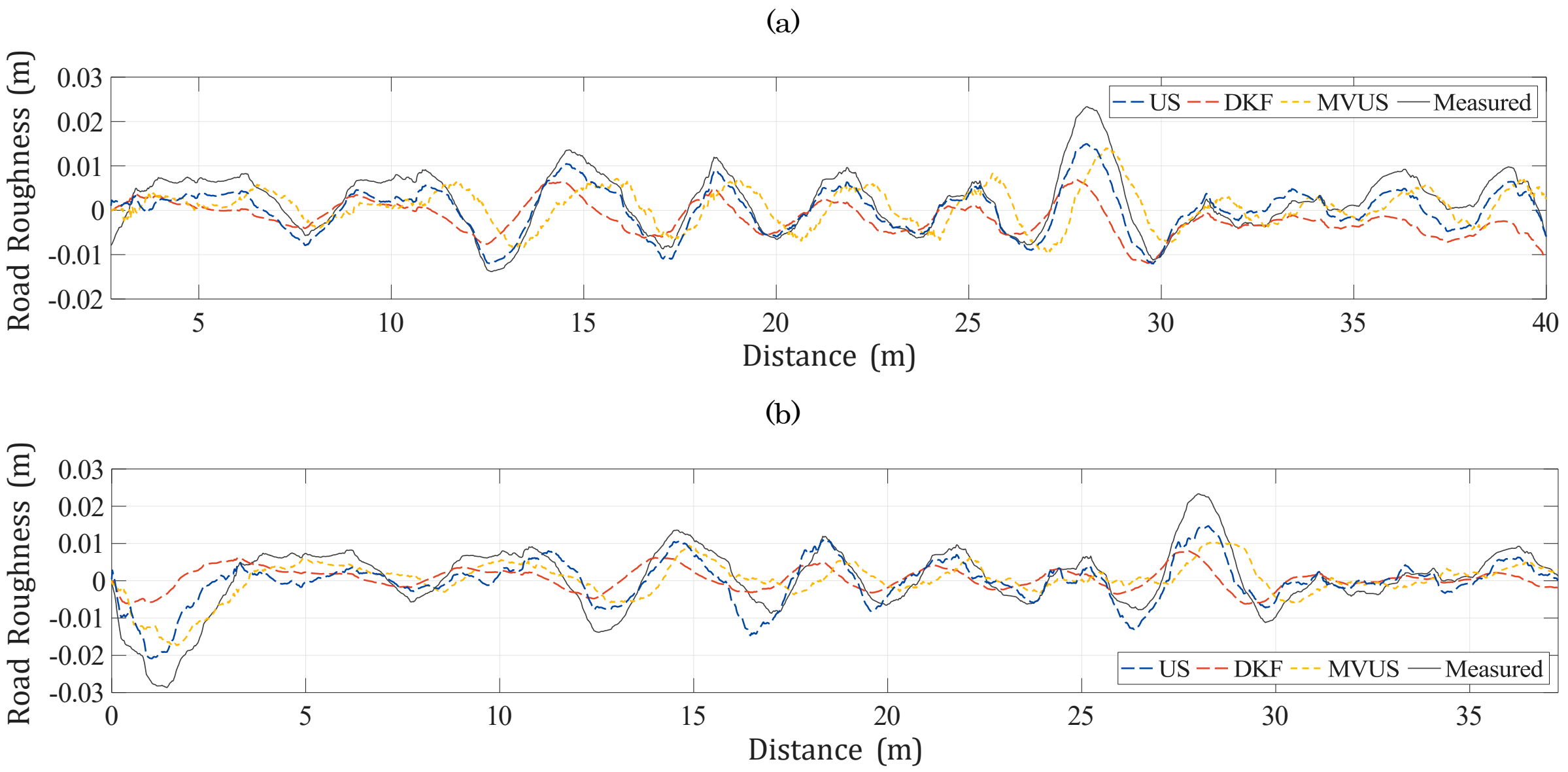

Figure 7: Identified Track 1 roughness profile with vehicle speed at approximately 10 km/h (Scenario 1): (a) front wheel, and (b) rear wheel.

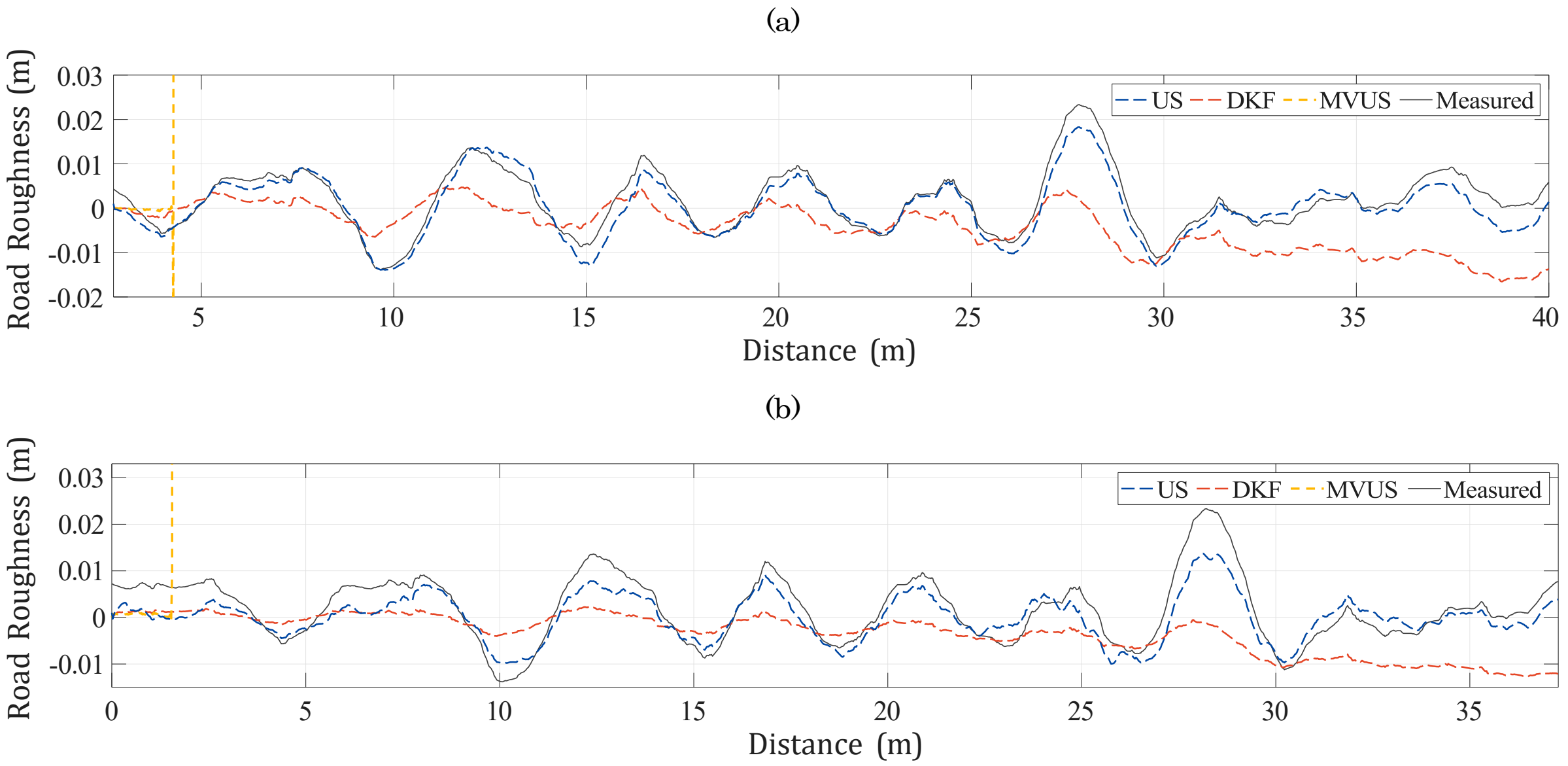

Figure 8: Identified Track 1 roughness profile with vehicle speed at approximately 20 km/h (Scenario 2): (a) front wheel, and (b) rear wheel.

The estimator behaviour is further illustrated through the propagation of input error covariance. Fig. 9a shows that the MVUS diverges within less than one second, leading to failure. In contrast, Fig. 9b demonstrates that the DKF covariance rapidly stabilises at a constant level, although small biases gradually accumulate and cause drift over time. Meanwhile, Fig. 9c indicates that the US converges to a bounded level with mild periodic oscillations, implying that accumulated integration errors are regularly corrected as the sliding window progresses. These observations underscore the advantage of a smoothing strategy under noisy conditions and highlight the robustness of the US against ill-conditioning.

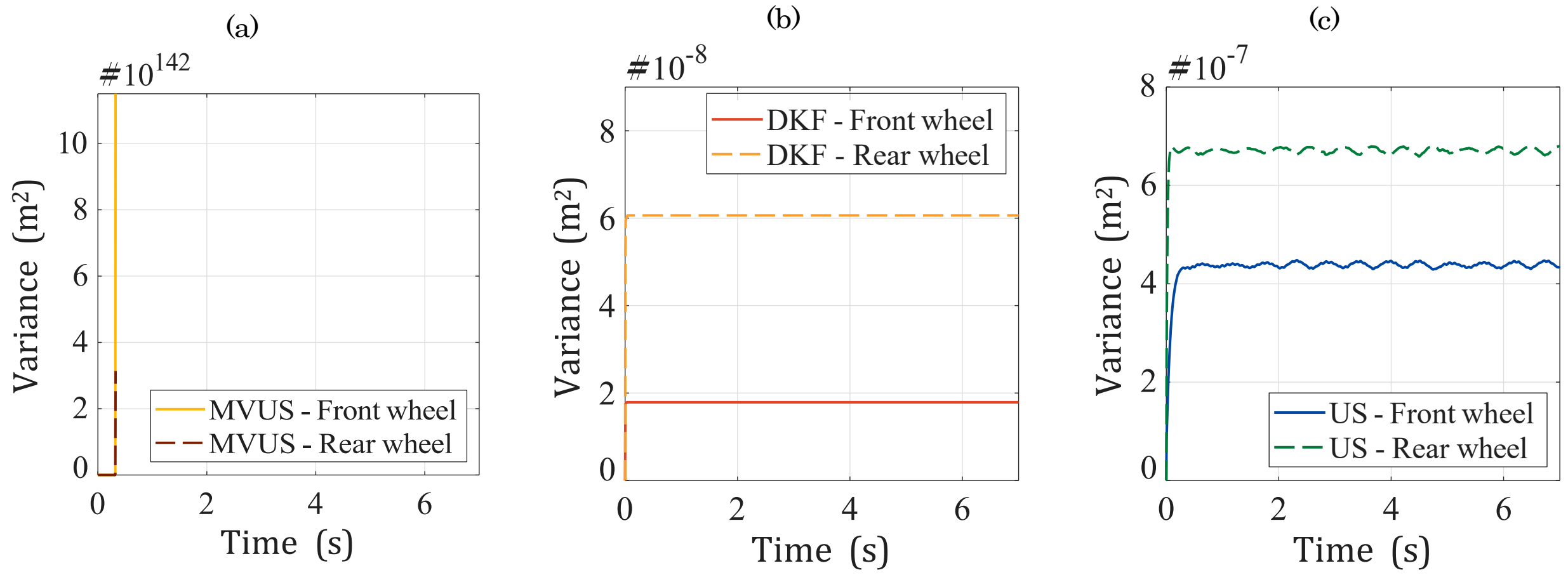

Figure 9: Input error covariance propagation during Scenario 2: (a) MVUS, (b) DKF, and (c) US.

Fig. 10 presents the Track 2 roughness profile identified with vehicle speed at approximately 10 km/h. All three estimators capture the overall profile shape, with the US most closely matching the measured profile across the test section and showing less underestimation at the local maxima around 17 m and 22 m. The DKF tends to underestimate the magnitude of sharper segments – more noticeably for the rear wheel – while the MVUS exhibits a small delay.

The spatial-frequency comparison is presented in Fig. 11, in which all estimators show a mismatch within the frequency range lower than $10^{-1}$ cycle/m. Over the mid-to-high spatial-frequency band, the US aligns with the measured magnitude and preserves more high-frequency content, whereas both DKF and MVUS display discrepancies at higher spatial frequencies. The corresponding NRMSE values are 0.215 for the US, 0.243 for the MVUS, and 0.309 for the DKF.

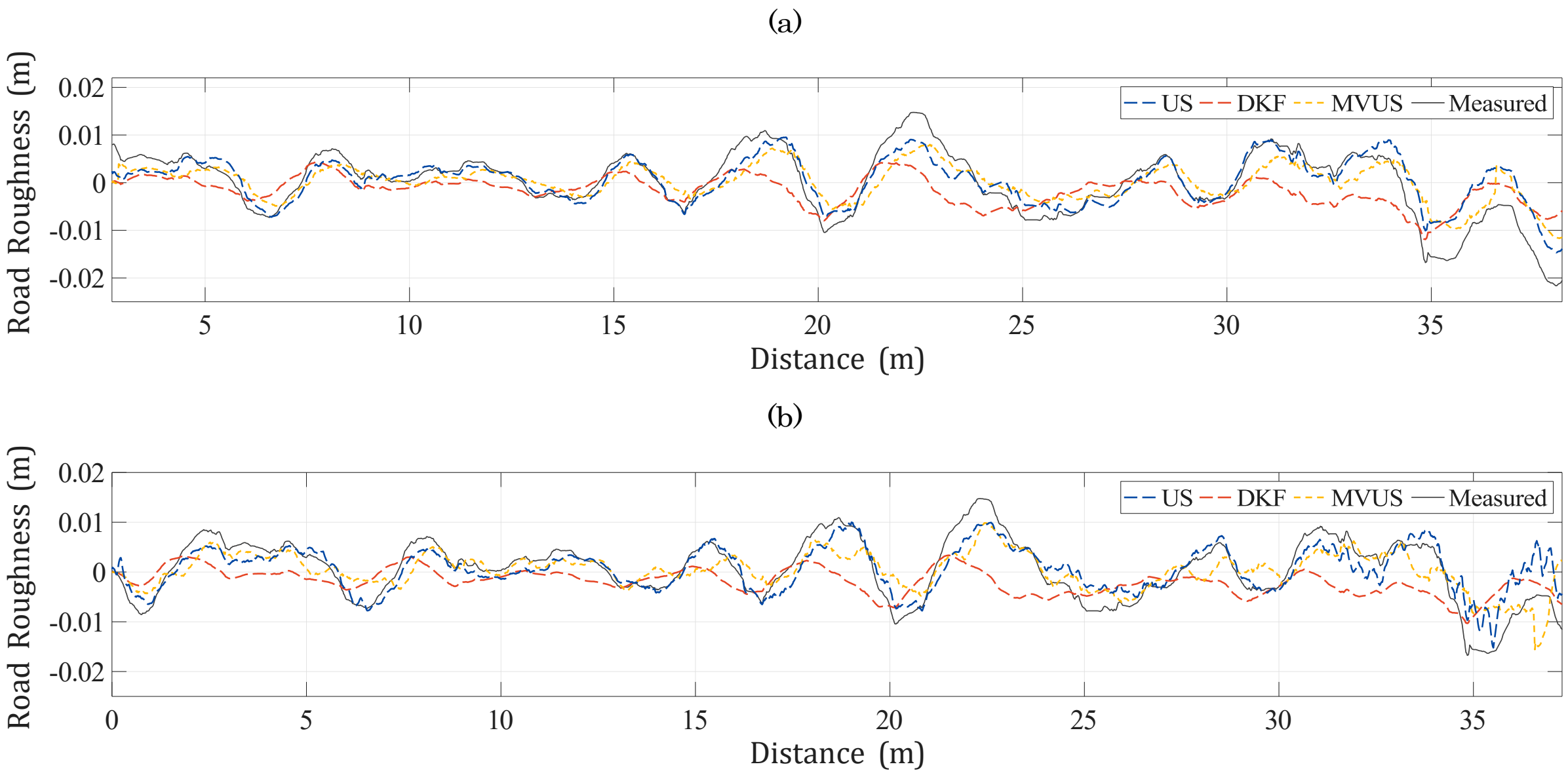

Figure 10: Identified Track 2 roughness profile with vehicle speed at approximately 10 km/h (Scenario 3): (a) front wheel, and (b) rear wheel.

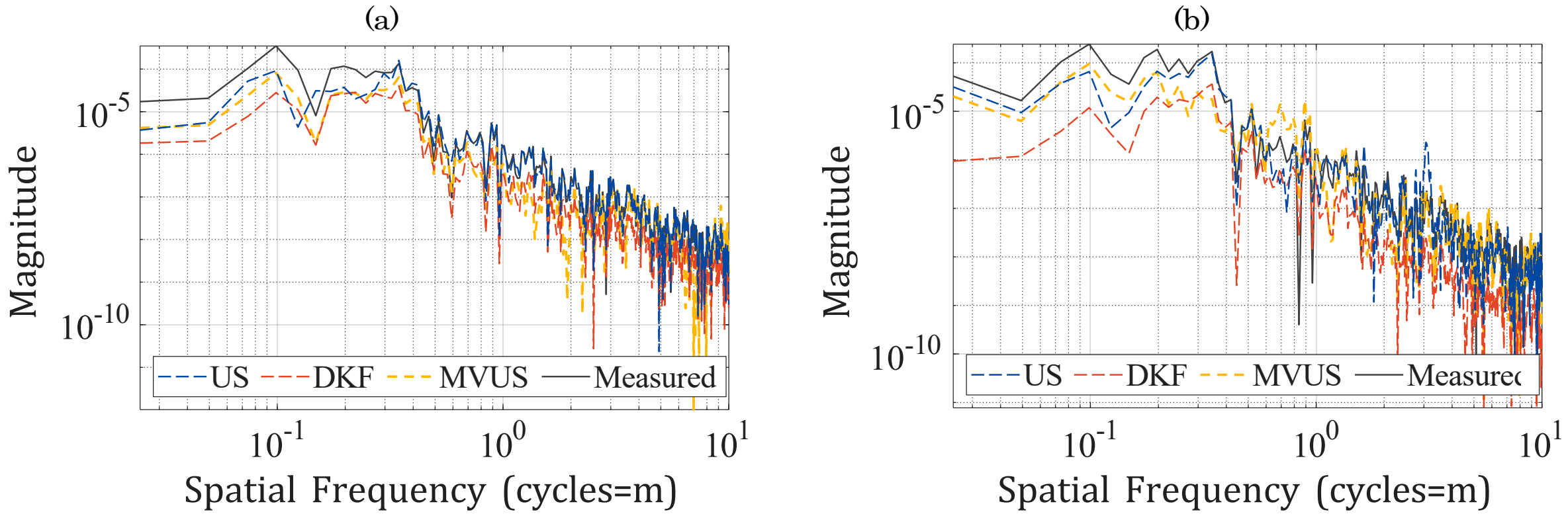

Figure 11: Spatial frequency of the identified Track 2 roughness profile with vehicle speed at approximately 10 km/h (Scenario 3): (a) front wheel, and (b) rear wheel.

Fig. 12 presents the Track 2 roughness profile identified with vehicle speed at approximately 20 km/h. The MVUS is not plotted since it encounters numerical errors, as in Scenario 2. The US follows the measured profile along the span and does not show a drift issue. In contrast, the DKF develops a drift issue: its estimate shifts downward with distance, as clearly visible in the traces under both wheels after approximately 25 meters. The spatial-frequency comparison in Fig. 13 reflects these behaviours: the drift in DKF is more obvious in the rear wheel results, as shown in Fig.12b. However, the US agrees with the measured magnitude over the mid-to-high frequency range. Similar to Scenarios 1 and 2, both estimators exhibit inconsistencies in the frequency range below $10^{-1}$ cycles/m. The same phenomenon was also reported for other system inversion methods used for road roughness identification [33]. The NRMSE for the US is 0.227.

Across the four scenarios, the US provided stable and accurate reconstructions using axle accelerations only. It remained well-conditioned when the MVUS encountered numerical errors (Scenarios 2 and 4) and did not exhibit the drift issue observed for the DKF at 20 km/h. Quantitatively, the US NRMSEs were 0.216, 0.218, 0.215, and 0.227 (Scenarios 1–4), yielding a variation of 2.2%. These results indicate that the US performance is consistent across tracks and speeds and robust to increased noise.

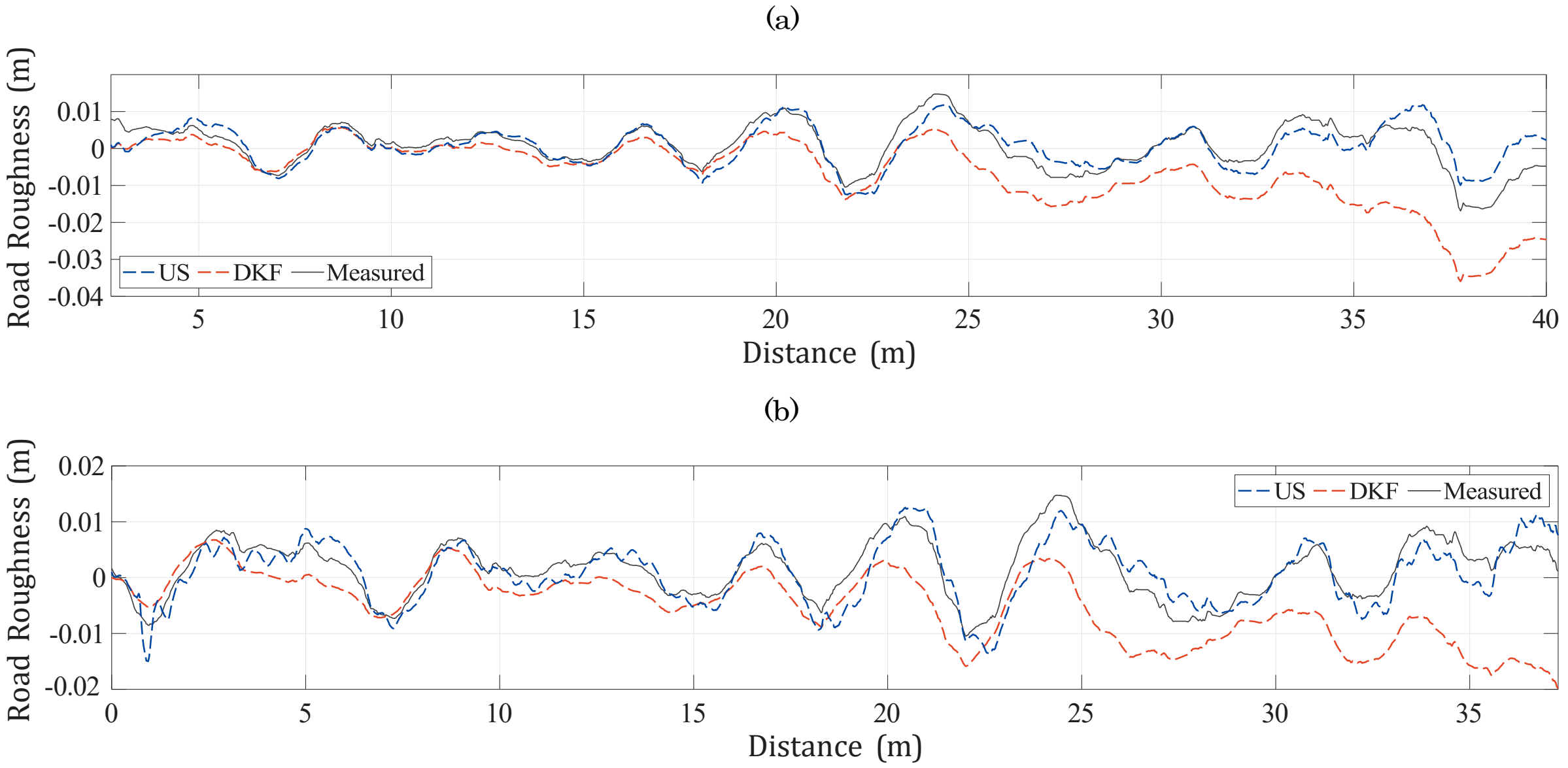

Figure 12: Identified Track 2 roughness profile with vehicle speed at approximately 20 km/h (Scenario 4): (a) front wheel, and (b) rear wheel.

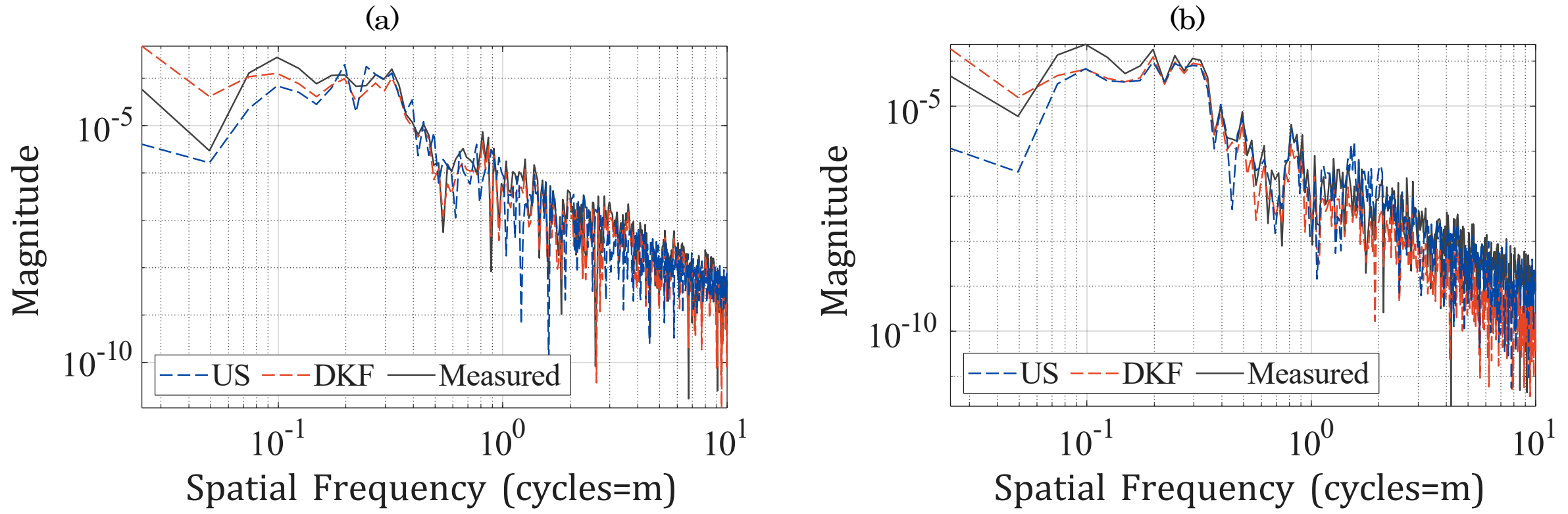

Figure 13: Spatial frequency of the identified Track 2 roughness profile with vehicle speed at approximately 20 km/h (Scenario 4): (a) front wheel, and (b) rear wheel.

## 5. Conclusions

This study introduced an output-only, drive-by road-roughness identification method – the Universal Smoothing (US) approach – using axle accelerations from a moving vehicle. A discrete half-car model with a backward-difference scheme brings both the current and previous roughness inputs into the observation equation, satisfying the inversion requirements of MVU-based estimators and enabling joint input–state estimation from limited measurements. A modified weight matrix was incorporated in the US weighted least squares to account for uncertainty in the previous-step input estimate, improving robustness when past roughness enters the current observation. Numerical stabilisation employed truncated singular value decomposition, and a practical grid-search strategy was used for tuning. Field data from a full-scale bridge test conducted at two speeds validated the approach against the Dual Kalman Filter (DKF) and the Minimum-Variance Unbiased Smoothing (MVUS).

- Across all scenarios, the US produced stable, accurate roughness reconstructions. Under higher measurement-noise conditions, the US remained well-conditioned, whereas the DKF exhibited drift and the MVUS encountered ill-conditioning – consistent with the identification plots and input estimation error covariance evolutions.

- A grid search tuned process noise and the number of retained singular values in the pseudoinverse. Truncated singular value decomposition governs the bias–variance trade-off: retaining too few singular values preserves stability but limits observability, whereas including very small singular values can lead to ill-conditioning. Accuracy gains diminish beyond a certain observation window length while computational cost grows rapidly.

- The method requires only axle-acceleration measurements and a simplified half-car model – no displacement instrumentation or vehicle modifications – supporting practical drive-by deployment.

Overall, the work establishes a solid foundation for reconstructing inputs at the tyre–road contact, paving the way toward more robust drive-by bridge health monitoring. Future work will incorporate explicit vehicle–bridge interaction to enable drive-by system identification of bridges while maintaining minimal instrumentation on the test vehicle.

## Acknowledgment

The authors thank Infrastructure Innovation Engineering at Kyoto University for providing the field data. This work was supported by the Australian Research Council (ARC) through a Discovery Early Career Researcher Award (DECRA, DE210101625) and the Industrial Transformation Research Program (IH210100048).